\documentclass{amsart}

\usepackage{amssymb,amsfonts,latexsym}

\numberwithin{equation}{section}

\newtheorem{theorem}{Theorem}[section]
\newtheorem{Theorem}[theorem]{Theorem}

\newtheorem{proposition}[theorem]{Proposition}

\newtheorem{Corollary}[theorem]{Corollary}

\theoremstyle{definition}

\newtheorem{definition}[theorem]{Definition}

\newtheorem{question}[theorem]{Question}

\newtheorem *{Theorem A}{Theorem A}
\newtheorem *{Theorem B}{Theorem B}
\newtheorem *{Corollary C}{Corollary C}

\newcommand{\bR}{{\bar R}}

\newcommand{\la}{\langle\,}

\newcommand{\ra}{\,\rangle}

\newcommand{\ben}{\begin{enumerate}}

\newcommand{\een}{\end{enumerate}}

\newcommand{\Jac}{{\rm Jac}}
\newcommand{\hT}{\widehat T}
\newcommand{\hR}{\widehat R}
\newcommand{\hS}{\widehat S}
\newcommand{\bS}{\bar S}
\newcommand{\bg}{\bar g}
\newcommand{\brf}{\bar f}
\newcommand{\im}{{\rm im}}
\newcommand{\res}{{\rm res}}
\newcommand{\Span}{{\rm Span}}

\begin{document}

\title
{Maximal Crossed Product Orders over Discrete Valuation Rings}

\author[Yuval Ginosar]
{Yuval Ginosar}

\date{\today}

\maketitle
\centerline {\small Department of Mathematics, University of Haifa,}
\centerline {\small Haifa 31905, Israel}
\centerline {\small \tt ginosar@math.haifa.ac.il}

\begin{abstract}
The problem of determining when a (classical) crossed product $T=S^f*G$
of a finite group $G$ over a discrete valuation ring $S$ is a
maximal order, was answered in the 1960's for the case where $S$ is
tamely ramified over the subring of invariants $S^G$. The answer was
given in terms of the conductor subgroup (with respect to $f$) of the inertia. In this
paper we solve this problem in general when $S/S^G$ is residually
separable. We show that the maximal order property entails a
restrictive structure on the sub-crossed product graded by the
inertia subgroup. In particular, the inertia is abelian. Using
this structure, one is able to extend the notion of the conductor.
As in the tame case, the order of the conductor is equal to the
number of maximal two sided ideals of $T$ and hence to the number of maximal orders
containing $T$ in its quotient ring. Consequently, $T$ is a maximal order
if and only if the conductor subgroup is trivial.
\end{abstract}

\section{Introduction}
Let $S$ be a discrete valuation ring (DVR) and let $G$ be a finite subgroup of Aut$(S)$.
Denote the unique maximal ideal of $S$ by $M_S$ and the corresponding residue field $S/M_S$ by $\bS$.
For any $f\in Z^2(G,S^*)$
consider the crossed product $T:=S^{f}*G=\oplus_{g\in G}SU_g$ with multiplication
\begin{equation}
sU_gtU_h=sg(t)f(g,h)U_{gh} \ \ s,t\in S,\ \ g,h\in G.
\end{equation}
Let $R:=S^G$ be the subring of $G$-invariant elements in $S$ and let $\bR:=R/(M_S\cap R)$ be its residue field.
We shall always assume that the extension $\bS/\bR$ is separable ({\it residual separability} property of $S/R$).
Denote the field of quotients of $S$ by $L$. Then the 2-cocycle $f$ can be regarded also as in $Z^2(G,L^*)$,
and $T$ is an $R$-order in the central simple algebra $L^{f}*G$.

\begin{question}\label{Q}
When is the $R$-order $T$ maximal in $L^{f}*G?$
\end{question}

Suppose that $S/R$ is {\it tamely ramified}, that is when
the order of the inertia subgroup $G_I\lhd G$ is prime to $p:=$char$(\bR)$.
In this case the answer to Question \ref{Q} can be given in terms of
the subgroup $H_f\lhd G$,
which is maximal in the inertia subgroup such that the cohomology class
$[\brf]\in H^2(G,\bS^*)$ is
inflated from $H^2(G/H_f,\bS^*)$, namely the {\it conductor} subgroup with respect to $f$.
\begin{Theorem}\cite[Theorem 2.5]{W63}
Let $S/R$ be a tamely ramified extension. Then the number of maximal $R$-orders containing
$T$ in $L^{f}*G$ is equal to the order of the conductor $H_f$.
In particular, $T$ is a maximal $R$-order if and only if $H_f$ is trivial.
\end{Theorem}
Question \ref{Q} is discussed in \cite{J} in a special instance of the tamely ramified case,
namely where $L$ is a finite extension of the $p$-adic rationals $\mathbb{Q}_p$.
The number of maximal $R$-orders containing $T$ in $L^f*G$
is given there in terms of the Schur index of the class $[f]\in H^2(G,L^*)$.

These results are generalized in \cite{CTA, TA} for any extension $S/R$ such that the
residue fields are finite. However, under this condition on the residue fields, $T$ cannot
be a maximal $R$-order unless $S/R$ is again tamely ramified (by \cite[Theorem 2]{H64}
and Theorem \ref{HAG} hereafter).\\

In this note we answer Question \ref{Q} dropping the above tame
ramification assumption. We first show
\begin{Theorem A}
{\it If $S^f*G$ is a maximal
order or, more generally, a hereditary $R$-order, then the inertia
subgroup of $G$ is abelian.}
\end{Theorem A}
With the restrictive structure that the heredity property entails on $T$ (Corollary \ref{ssiff}),
we are able to extend the notion of the conductor subgroup
(Definition \ref{condef}).
This notion arises naturally from a well known group cohomology map.
It turns out that the image of this map controls the number of maximal two sided ideals of $T$
or, equivalently, the number of maximal orders in $L^f*G$ which contain $T$
(Corollary \ref{image}).
This implies that as in the tamely ramified case, the maximal order property of $T$ depends on the
triviality of the conductor:
\begin{Theorem B}
{\it Let $T=S^f*G$ be a hereditary crossed product order. Then the number of maximal $R$-orders containing
$T$ in $L^{f}*G$ is equal to the order of the conductor $H_f$.
In particular, $T$ is a maximal $R$-order if and only if
it is hereditary and the conductor $H_f$ is trivial.}
\end{Theorem B}

Finally, the demand that $S$ is a DVR can be relaxed to the more general case where $S$ is a Dedekind domain,
$R=S^G$ is a DVR and $S/R$ is residually separable. The reduction is fairly standard and appears in Section \ref{reduct}.\\

\noindent{\bf Acknowledgement.} I am grateful to my colleague Amiram Braun and to
Mota Frances for their help.

\section{Heredity and Semisimplicity}\label{hs}
The following result will be useful in the sequel.
\begin{Theorem}\label{HAG}\cite[Theorem 2.3]{AG}
Let $R$ be a DVR and let $\Lambda$ be an $R$-order.
Then $\Lambda$ is maximal if and only if it is hereditary
and has a unique two sided ideal.
\end{Theorem}


We first handle the heredity condition in Theorem \ref{HAG}.
In order to formulate the criterion, note that $T/M_S T$ is isomorphic to the crossed product
$\bS^{\brf}*G$. The action of $G$ on $\bS$ is induced by its action on $S$ and hence admits a
kernel. This kernel is the inertia (or the first ramification) subgroup $G_I$. The 2-cocycle
$\brf$ is the image of $f$ under the natural map $Z^2(G,S^*)\to Z^2(G,\bS^*)$. We have
\begin{Theorem}\label{heredity}\cite[Theorem A]{BGL}
With the above notation,
$S^f*G$ is a hereditary order if and only if $\bS^{\brf}*G$ is semi-simple.
\end{Theorem}

If $S/R$ is tamely ramified, then the fact that the order of $G_I$
is invertible in the field $\bS$ implies that $\bS^{\brf}*G$ is semi-simple independently of $\brf$,
by a generalized Maschke's Theorem.
Hence $T$ is hereditary.
However, it turns out that $T$ may be hereditary even when $S/R$ is not tamely ramified
\cite[Example 4.1]{BGL}.

Here is an explicit criterion for the semisimplicity of $\bS^{\brf}*G$.
By Theorem \ref{heredity},
it is a necessary and sufficient condition for the heredity property of $S^{f}*G$.
Note that since the inertia subgroup $G_I\lhd G$ acts trivially on $\bS$, the sub-crossed product graded by $G_I$
is a twisted group algebra ${\bS}^{\brf}G_I$.
Let $P$ be a $p$-Sylow subgroup of $G_I$,
where $p$ is the characteristic of the residue field $\bS$
(in case $p=0$, take $P$ as the trivial group).
\begin{Theorem}\label{AR}\cite[Theorem 2]{AR},
With the above notation, $\bS^{\brf}*G$ is semi-simple
if and only if the twisted group subalgebra $F:=\bS^{\brf}P$ is a purely inseparable field extension of $\bS$.
In particular, $P$ is abelian and the 2-cocycle $\brf$ is non-trivial on any non-trivial subgroup of $P$.
Additionally, it follows that the order of the commutator subgroup $[G_I,G_I]$ is prime to $p$.
\end{Theorem}

\noindent {\it Proof of Theorem A}.
Let $w$ be a generator of $M_S$. For any $\sigma\in G_I$ let $\sigma(w)=x_{\sigma}w$,
where $x_{\sigma}\in S^*$.
Then by \cite[Theorem 25, P. 295]{ZS}, the map $\sigma\mapsto {\bar x_\sigma}$ is a homomorphism from $G_I$ into $\bS^*$
whose kernel under the residual separability assumption is exactly $P$ (the second ramification group).
Consequently, $P$ is normal and
$G_I=P\rtimes C_{e_0}$, where $C_{e_0}=\la \sigma_0\ra $
is a cyclic group whose order
is prime to $p$. Now, if $\bS^{\brf}*G$ is semisimple, then by Theorem \ref{AR}, the order
of the commutator $[G_I,G_I]$ is prime to $p$, and hence the action of $C_{e_0}$ on $P$ is
trivial. Consequently, $G_I$ is a direct product of $P$ and $C_{e_0}$, hence abelian.
\qed

The following is a stronger consequence of the semisimplicity of $\bS^{\brf}*G$.
For the sake of convenience, we continue to denote the basis elements of $\bS^{\brf}*G$ by $\{U_{\sigma}\}_{\sigma\in G}$.
For any $\sigma\in P$, let
$\lambda:=U_{\sigma}U_{\sigma_0}U_{\sigma}^{-1}U_{\sigma_0}^{-1}\in
\bS^*$. Suppose that the order of $\sigma$ is $p^m$ for some $m$.
Then
$\lambda^{p^m}:=U_{\sigma}^{p^m}U_{\sigma_0}U_{\sigma}^{-p^m}U_{\sigma_0}^{-1}=1$.
Since $\bS$ does not admit non-trivial $p$-th roots of 1, we deduce
that $\lambda=1$ and thus $U_{\sigma_0}$ is central in
${\bS}^{\brf}G_I$.
Let $\alpha_0:=U_{\sigma_0}^{e_0}\in \bS^*$. We obtain
\begin{Corollary}\label{ssiff}
$T$ is hereditary if and only if ${\bS}^{\brf}G_I$ is semisimple
and isomorphic to a commutative
twisted group algebra $F^{\alpha_0}C_{e_0}\simeq F[x]/\la x^{e_0}-\alpha_0\ra $,
where $F=\bS^{\brf}P$ is a purely inseparable extension of
$\bS$ and $C_{e_0}$ is a cyclic group of order $e_0$, which is prime to $p$.
In particular, $G_I$ is abelian of the form $G_I=P\times C_{e_0}$.
\end{Corollary}

\section{The Number of Simple Components of $\bS^{\brf}*G$}

Suppose that $T$ is hereditary. Then by Theorems \ref{heredity} and \ref{AR},
the restriction of $\brf$ to a subgroup $H$ of $G_I$
can be trivial only if $H$ is a $p'$-group,
that is $H$ is contained in $C_{e_0}$ and hence normal in $G$.
Due to this observation, one can generalize the
definition of the conductor subgroup as follows.
\begin{definition}\label{condef}
Let $f\in Z^2(G,S^*)$ such that $S^f*G$ is hereditary.
The conductor $H_f$ with respect to $f$ is the maximal subgroup of the inertia
such that the class $[\brf]\in H^2(G,\bS^*)$ is inflated from $H^2(G/H_f,\bS^*)$.
\end{definition}
In subsection \ref{lines} we make use of the above definition so as to
obtain that the number of simple components of
$\bS^{\brf}*G$ and the number of maximal two-sided ideals in $T$
are both equal to the order of the conductor subgroup $H_f$
(compare with \cite[Theorem 2.5]{W63}).
By that Theorem B will be deduced,
since the number of maximal $R$-orders containing $T$ in $L^f*G$
is equal to the number of maximal two-sided ideals in $T$ \cite[Theorem 1.7]{H63}.

The proof of Theorem B is partially based on \cite{W63}.
Subsection \ref{coho} below proposes a cohomological
interpretation to this result.

\subsection{}\label{coho}
In this subsection we present the cohomological tool for the calculation of the
number of simple components of $\bS^{\brf}*G$ that is essential for Theorem B.
The discussion is based on a construction due to J.P. Serre
and can be found in \cite[Section 1.7]{K}. Here is a brief description.

Let
\begin{equation}\label{ext}
1\rightarrow A \rightarrow G \rightarrow G/A \rightarrow 1
\end{equation}
be an extension of finite groups, where $A$ is abelian.
As usual, $G/A$ acts on $A$ via the conjugation in $G$, namely, for
every $\bar{g}\in G/A$ and $a\in A$, $\bg(a)=g ag^{-1}$. Next, let $M$
be a $G$-module which is $A$-trivial, that is $G$ acts on $M$ via
$G/A$. Then the action of $G/A$ on $A$ induces the following diagonal action of
$G/A$ on $\hom(A,M)(\simeq H^1(A,M))$.
Let $\bg\in G/A$ and $\varphi\in \hom(A,M)$. Then
$\bg(\varphi)\in \hom(A,M)$ is defined on $a\in A$ via the pairing
\begin{equation}
\la \bg(\varphi),a\ra =\bg\la \varphi,\bg^{-1}(a)\ra .
\end{equation}

Next, let $f\in Z^2(G,M)$ satisfy
\begin{equation}\label{karp}
f(g_1,g_2)=f(g_1,g_2a),\ \ \forall g_1,g_2\in G, a\in A.
\end{equation}
In particular, the restriction of $f$ to $A$ is trivial.

For any $a\in A$ and $\bg\in G/A$ define
\begin{eqnarray}\label{pi}
\begin{array}{rcl}
\pi_f(\bg):A & \rightarrow & M\\
a & \mapsto & f(a,g).
\end{array}
\end{eqnarray}
\begin{Theorem}\label{items}(see \cite[Theorem 7.3, P. 60]{K})
Let $\res^G_A:H^2(G,M)\to H^2(A,M)$ and $\inf^{G/A}_G:H^2(G/A,M)\to H^2(G,M)$
be the restriction and inflation maps respectively.
Let $\pi_f$ be as in (\ref{pi}). Then
\begin{enumerate}
\item Any class in $\ker(\res^G_A)$ admits a representative $f\in Z^2(G,M)$ satisfying (\ref{karp}).
\item For any $\bg\in G/A$, $\pi_f(\bg)\in \hom(A,M)$.
\item $\pi_f(\bg)$ does not depend on the choice of the representative $g\in G$ for $\bg$.
\item The map $\bg\mapsto\pi_f(\bg)$ is a 1-cocycle from $G/A$ to $\hom(A,M)$.
\item If $f'\in[f]$ satisfies (\ref{karp}), then the 1-cocycles $\pi_{f'}$ and $\pi_f$ differ by a 1-coboundary.
\item If $f_1$ and $f_2$ satisfy (\ref{karp}), then so does $f_1+f_2$. Moreover, $\pi_{f_1+f_2}=\pi_{f_1}+\pi_{f_2}.$
\item $\pi_f\in B^1(G/A,\hom(A,M))$ if and only if the cohomology class $[f]$ is in the image of $\inf^{G/A}_G$.
\end{enumerate}
\end{Theorem}
\begin{Corollary}(see \cite[Theorem 7.3, P. 60]{K})
The map $\Pi:[f]\mod[\im(\inf^{G/A}_G)]\mapsto[\pi_f]$ is a well defined injection of
$\ker(\res^G_A)/\im(\inf^{G/A}_G)$ into $H^1(G/A,\hom(A,M))$.
\end{Corollary}

The map $\Pi$ is applied for crossed products as follows.
Let $K^f*G=\oplus_{g\in G}KU_g$ be a crossed product, where $K$ is a field and $f\in Z^2(G,K^*)$.
Suppose that $A\vartriangleleft G$ is an abelian subgroup acting trivially on $K$
such that the restriction of $f$ to $A$ is cohomologically trivial.
By Theorem \ref{items}(1), the $K$-basis $\{U_g\}_{g\in G}$ may be chosen such that
\begin{equation}\label{3.5}
U_{ga}=U_gU_a, \forall g\in G, a\in A.
\end{equation}
In particular, $K^f*G$ contains the ordinary group algebra $KA$.
Then $G/A$ acts on $KA$ via the conjugation in $K^f*G$.
We describe this action using the 1-cocycle $\pi_f\in Z^1(G/A,\hom(A,K^*))$.
Let $ \bg\in G/A$ and $a\in A$. Then by (\ref{pi}),
\begin{equation}
\la \pi_f(\bg),a\ra =f(a,g)=U_aU_{g}U_{ag}^{-1}=U_aU_{g}U_{g \bg^{-1}(a)}^{-1}.
\end{equation}
By (\ref{3.5}), $U_{g \bg^{-1}(a)}^{-1}=U_{\bg^{-1}(a)}^{-1}U_{g}^{-1}$. Consequently,
\begin{equation}\label{3.7}
\la \pi_f(\bg),a\ra =U_aU_{g}U_{\bg^{-1}(a)}^{-1}U_{g}^{-1}.
\end{equation}
Hence, for every $ \bg\in G/A$ and $a\in A$
\begin{equation}\label{act}
\bg(U_a)=U_{g}U_aU_{g}^{-1}=(U_{g}U_{\bg^{-1}(\bg(a))}^{-1}U_{g}^{-1})^{-1}=\la \pi_f(\bg),\bg(a)\ra ^{-1}U_{\bg(a)}.
\end{equation}
Now, suppose that $|A|$ is invertible in $K$. Then the primitive idempotents of $KA$ are
$\iota_{\chi}=\frac{1}{|A|}\sum_{a\in A}\la \chi,a\ra ^{-1}U_a$ for every $\chi\in \hom(A,K^*)$.
The action on $KA$ yields an action of $G/A$ on the set of primitive idempotents of $KA$ as follows.
\begin{proposition}\label{BG}(see a special instance in \cite[Proposition 2.9]{BG})
With the above notation, let $\bg\in G/A$ and let $\chi\in\hom(A,K^*)$.
Then $\bg(\iota_{\chi})=\iota_{\bg(\chi)\pi_f(\bg)}$.
\end{proposition}
\begin{proof}
\begin{equation*}
\begin{split}
\bg(\iota_{\chi})=U_{g}\iota_{\chi}U_{g}^{-1}=U_{g}\frac{1}{|A|}\sum_{a\in A}\la \chi,a\ra ^{-1}U_aU_{g}^{-1}
=\frac{1}{|A|}\sum_{a\in A}\bg\la \chi,a\ra ^{-1}U_{g}U_aU_{g}^{-1}.
\end{split}
\end{equation*}
Then by (\ref{act}),
\begin{equation*}
\begin{split}
\bg(\iota_{\chi})~=&~\frac{1}{|A|}\sum_{a\in A}\bg\la \chi,a\ra ^{-1}\la \pi_f(\bg),\bg(a)\ra ^{-1}U_{\bg(a)}\\
=&~\frac{1}{|A|}\sum_{a\in A}\la \bg(\chi),\bg(a)\ra ^{-1}\la \pi_f(\bg),\bg(a)\ra ^{-1}U_{\bg(a)}\\
=&~\frac{1}{|A|}\sum_{a\in A}\la \bg(\chi)\pi_f(\bg),a\ra ^{-1}U_{a}=\iota_{\bg(\chi)\pi_f(\bg)}.
\end{split}
\end{equation*}
\end{proof}

\subsection{}\label{lines}
The second step in determining if $T$ is a maximal order, after having taken care of its heredity
property (in Section \ref{hs}), is to handle the locality condition in Theorem \ref{HAG}.
We have
\begin{proposition} \label{local}
The number of maximal two-sided ideals in $T$ is equal to the number of maximal two-sided ideals in $\bS^{\brf}*G$.
In particular, $T$ is local if and only if so is $\bS^{\brf}*G$.
\end{proposition}
\begin{proof}
This is clear since every maximal two sided ideal of $T$ contains $M_ST$.
\end{proof}
Assume that $\bS^{\brf}*G=\Span_{\bS}\{U_{g}\}_{g\in G}$ satisfies the necessary and sufficient
condition for semisimplicity in Corollary \ref{ssiff}.
Then by Proposition \ref{local}, the number of maximal two-sided ideals in $T$
is equal to the number of simple components of $\bS^{\brf}*G$.
In particular, by Theorem \ref{HAG}, $T$ is a maximal order if and only if $\bS^{\brf}*G$
admits a single simple component.

We need to deal with the following
\begin{question}\label{simples}
Let ${\bS}^{\brf}*G=T/M_ST$ be a crossed product as above.
Suppose that ${\bS}^{\brf}*G$ is semisimple. How many simple components does ${\bS}^{\brf}*G$ admit?
In particular, when is  ${\bS}^{\brf}*G$ simple?
\end{question}
In general, determining the number of simple components of an
arbitrary semisimple crossed product $K^f*G$ of a finite group $G$ over a
field $K$ might be hard. Suppose that $[f]\in\ker(\res^G_A)$
for an abelian subgroup $A\lhd G$ which acts trivially on $K$
(and by Theorem \ref{items}(1) we may assume that $f$ satisfies (\ref{karp})). Then
a necessary condition for the simplicity of $K^f*G$ is that the primitive idempotents of the
commutative group ring $KA$ belong to the same orbit under the action of $G$.
By Proposition \ref{BG}, this implies that the 1-cocycle $\pi_f$ is {\it onto} $\hom(A,K^*)$.
Under our conditions however, the central idempotents of ${\bS}^{\brf}*G$ can be calculated using Proposition \ref{BG},
as well as the structure of ${\bS}^{\brf}*G_I$ given in Corollary \ref{ssiff}.

The following claim shows that the central primitive idempotents of $\bS^{\brf}*G$
are supported by the inertia subgroup.
\begin{proposition}\label{H}
The center of ${\bS}^{\brf}*G$ lies in ${\bS}^{\brf}G_I$.
\end{proposition}
\begin{proof}
Let $y=\sum_{g\in G}{\bar s_g}U_g\in {\bS}^{\brf}*G$.
Suppose that ${\bar s}_{g_0}\neq 0$ for some $g_0 \notin G_I$.
Then since $g_0$ is not in the kernel of the action of $G$ on $\bS$,
there exists an element $\bar{t}\in{\bS}$ which
does not commute with $U_{g_0}$ and hence also with $y$.
\end{proof}

In view of Proposition \ref{H},
any central idempotent of $\bS^{\brf}*G$ is a sum of certain primitive idempotents of
the commutative twisted group subalgebra ${\bS}^{\brf}G_I$.
By Corollary \ref{ssiff}, ${\bS}^{\brf}G_I$
is isomorphic to the commutative twisted group ring
$F^{\alpha_0}C_{e_0}=\Span_{F}\{U_{\sigma_0^i}\}_{0\leq i\leq e_0-1}$,
where $(U_{\sigma_0})^{e_0}=\alpha_0\in \bS^*$.
We need the following properties of the field $\bS$.
\begin{proposition}\label{eroots}
With the above notation
\begin{enumerate}
\item The field $\bS$ contains all $e_0$-th roots of $1$.
\item Let $\zeta_{e_0} \in \bS$ be an $e_0$-th root of $1$.
Then for every $g\in G$, $g\sigma_0g^{-1}=\sigma_0^m$ where $m$ is determined by
$g(\zeta_{e_0})=\zeta_{e_0}^m$.
\end{enumerate}
\end{proposition}
\begin{proof}
The map $\sigma\mapsto {\bar x_\sigma}$ in the proof of Theorem A
yields an embedding of the cyclic group $C_{e_0}$ in $\bS^*$ verifying (1).
In order to prove (2), we need to show that this map is also a $G$-morphism.
As can easily be seen, the map does not depend on the generator $w$ of $M_S$.
Choosing $g^{-1}(w)$ as a new generator we obtain that $\sigma(g^{-1}(w))=y_{\sigma}g^{-1}(w)$, where
${\bar y_\sigma}={\bar x_\sigma}$. Acting with $g$ on both sides gives $g\sigma g^{-1}(w)=g(y_{\sigma})w$.
Hence ${\bar x_{g\sigma g^{-1}}}=g({\bar x_{\sigma}})$ and we are done.
\end{proof}

Let $\Gamma_f$ be a maximal subgroup of $G_I$
such that the restriction of $\brf$ to it is
cohomologically trivial (compare with \cite[P. 111, Definition]{W63}).
By Theorem \ref{AR},  $\Gamma_f$ intersects $P$ trivially, hence it is contained in $C_{e_0}=\la \sigma_0\ra $
and therefore it is unique.
Let $c=c(f)$ be such that $\Gamma_f:=\la \sigma_0^{c}\ra $.
Then the order of $\Gamma_f$ is $d=d(f)=\frac{e_0}{c}$, which is the maximal divisor of $e_0$ such that $\alpha_0$
admits a root of order $d$ in $\bS$ (equivalently in $F$, since $\bS$ is its separable closure inside $F$
and $(e_0,p)=1$). In particular, $d$ is invertible in $\bS$.

It is clear that every subgroup $A$ of $\Gamma_f$ is normal in $G$, since these subgroups
are contained in the cyclic normal subgroup $C_{e_0}$.
The map $\pi_{\brf}$ may therefore be applied for every $A\lhd \Gamma_f$ and $M:=\bS^*$.
By Theorem \ref{items}(1), putting $A:=\Gamma_f$,
we may assume that $\brf$ satisfies (\ref{karp}).
In particular, $U_{\sigma_0^{cj}}=U_{\sigma_0^{c}}^j$ for every integer $j$.
Now, by Proposition \ref{eroots}(1), $\bS$ contains a primitive $d$-th roots of $1$,
denoted by $\zeta_d$. Let $k$ be a divisor of $d$, and let
$A:=\la \sigma_0^{kc}\ra $ be a subgroup of $\Gamma_f$ of order $\frac{d}{k}$.
Then $\hom(A,M)=\hom(A,\bS^*)$ is a cyclic group of order $\frac{d}{k}$ whose elements
are determined by the generator ${\sigma_0}^{kc}$ as follows.
\begin{equation}\label{charac}
\chi^{(k)}_j:{\sigma_0}^{kc}\mapsto \zeta_d^{kj} , \ \ 0\leq j\leq \frac{d}{k}-1.
\end{equation}
The idempotents of ${\bS}^{\brf}G_I$ can now be given explicitly.
For $A=\Gamma_f$ put $k=1$ and let $\chi_j=\chi^{(1)}_j$ in (\ref{charac}).
\begin{proposition}\label{idem}
The elements
\begin{equation}\label{iota}
\iota_j=\frac{1}{d}\sum_{l=0}^{d-1}\la \chi_j,\sigma_0^{cl}\ra ^{-1}U_{\sigma_0^{cl}}, \ \ 0\leq j \leq d-1
\end{equation}
form a complete set of primitive orthogonal idempotents of ${\bS}^{\brf}G_I$.
\end{proposition}
\begin{proof}
Since ${\bS}^{\brf}G_I\simeq F^{\alpha_0}C_{e_0}$,
one can apply \cite[Proposition 2.2]{W63} putting $F$ as the base field.
\end{proof}
The number of simple components of ${\bS}^{\brf}G$ depends on the action of $G/\Gamma_f$ on the above idempotents.
We have
\begin{proposition}\label{actriv}
For every $A=\la \sigma_0^{kc}\ra \lhd \Gamma_f$, the action of $G/A$ on $\hom(A,\bS^*)$ is trivial.
\end{proposition}
\begin{proof}
Let $\bg\in G/A$ and suppose that $\bg^{-1}({\sigma_0})=g^{-1}{\sigma_0}g={\sigma_0}^{m}$.
Then by Proposition \ref{eroots}(2), $\bg^{-1}(\zeta^k_{d})=g^{-1}(\zeta^k_{d})=\zeta_{d}^{km}$.
Now, let $\chi^{(k)}_j\in \hom(A,\bS^*)$. Then
$\la \bg(\chi^{(k)}_j),{\sigma_0}^{kc}\ra =\bg\la \chi^{(k)}_j,\bg^{-1}({\sigma_0}^{kc})\ra =\bg\la \chi^{(k)}_j,{\sigma_0}^{kcm}\ra
=\bg(\zeta_d^{kmj})=\zeta_d^{kj}=\la \chi^{(k)}_j,{\sigma_0}^{kc}\ra $, proving that $\bg(\chi^{(k)}_j)=\chi^{(k)}_j$.
\end{proof}
By Propositions \ref{BG} and \ref{actriv} for $A=\Gamma_f$, we obtain that an element $\bg\in G/\Gamma_f$
acts on the idempotents of ${\bS}^{\brf}G_I$ as translations by $\pi_{\brf}(\bg)$. More precisely
\begin{equation}\label{translate}
\bg(\iota_{\chi_j})=\iota_{\chi_j\pi_{\brf}(\bg)},\ \  \bg\in G/\Gamma_f, \ \ 0\leq j\leq d-1.
\end{equation}
By Proposition \ref{actriv}, the 1-cocycle $\pi_{\brf}: G/\Gamma_f\to\hom(\Gamma_f,\bS^*)$ is in fact a group homomorphism.
By (\ref{translate}), there is a 1-1 correspondence between
the orbits induced by the action of $G/\Gamma_f$ on the set of primitive idempotents of ${\bS}^{\brf}G_I$
and the cosets of the image $\pi_{\brf}(G/\Gamma_f)$ in $\hom(\Gamma_f,\bS^*)$ (and hence all the orbits are
of the same cardinality).

Next, by Propositions \ref{H} and \ref{idem}, any central idempotent of $\bS^{\brf}*G$ is of the form
$\iota=\sum_{j\in B}\iota_j$, where $B\subset\{0,...,d-1\}$ is a set of indices of
an orbit of primitive idempotents of ${\bS}^{\brf}G_I$ under the action of $G/\Gamma_f$.

Here is an answer to Question \ref{simples} in terms of the image of $\pi_{\brf}$.
\begin{Corollary}\label{image}
Let $T=S^f*G$ be a hereditary crossed product.
Then the number of simple components of $T/M_ST=\bS^{\brf}*G$ is equal to the index of $\pi_{\brf}(G/\Gamma_f)$ in $\hom(\Gamma_f,\bS^*)$.
In particular, $\bS^{\brf}*G$ is simple if and only if $|\pi_{\brf}(G/\Gamma_f)|=|\hom(\Gamma_f,\bS^*)|=d$.
\end{Corollary}
We now show that the number of simple components of $\bS^{\brf}*G$,
which is the same as the number of maximal orders containing $T$ in $L^f*G$,
is equal to the order of the conductor.\\
\noindent {\it Proof of Theorem B}.
Let $I_f$ be the index of $\pi_{\brf}(G/\Gamma_f)$ in $\hom(\Gamma_f,\bS^*)$ as above.
Then $I_f$ is the order of the maximal subgroup $A\lhd \Gamma_f$ such that the decomposition
\begin{equation}
G/\Gamma_f\xrightarrow{\pi_{\brf}}\hom(\Gamma_f,\bS^*)\xrightarrow{\res}\hom(A,\bS^*)
\end{equation}
is trivial, where the right map is the restriction map from $\Gamma_f$ to $A$.
Now, consider the diagram
\begin{eqnarray}\label{diag}
\begin{array}{ccc}
G/A & \longrightarrow & G/\Gamma_f\\
\downarrow & & \downarrow \\
\hom(A,\bS^*) & \xleftarrow{\res  } & \hom(\Gamma_f,\bS^*)
\end{array},
\end{eqnarray}
where the vertical arrows stand for the maps $\pi_{\brf}$ with respect to the normal subgroups $A$ and $\Gamma_f$
and the upper horizontal map is the natural projection.
By the definition of $\pi_{\brf}$ (\ref{pi}) and Theorem \ref{items}(3), we deduce that the diagram (\ref{diag}) is commutative.
Consequently, $A$ is the maximal subgroup of $\Gamma_f$ such that the map $G/A \xrightarrow{\pi_{\brf}} \hom(A,\bS^*)$
is trivial. Equivalently, since the action of $G/A$ on $\hom(A,\bS^*)$ is trivial (Proposition \ref{actriv}),
$A$ is the maximal subgroup of $\Gamma_f$ such that $\pi_{\brf}$ with respect to $A$ is a 1-coboundary.
By Definition \ref{condef} and Theorem \ref{items}(7), we obtain that $A$ coincides with the conductor $H_f$.
Applying Corollary \ref{image}, we obtain that the number of simple components of $\bS^{\brf}*G$ is equal to the
order of $H_f$. By Proposition \ref{local}, $|H_f|$ is the number of maximal two-sided ideals in $T$
and by \cite[Theorem 1.7]{H63}, it is the number of maximal $R$-orders containing $T$ in $L^f*G$.
This completes the proof of Theorem B.
\qed

\section{A Reduction Argument}\label{reduct}
In this section we reduce the question of when $T=S^f*G$ is a maximal order from the case where $S$ is a Dedekind
domain and $R$ is a DVR (keeping the demand that $S/R$ is residually separable)
to the case discussed in the previous sections, namely where both $S$ and $R$ are DVR's.

Let $z$ generate the unique maximal ideal of $R$ and let $\hR:=\lim\limits_{\leftarrow i} R/z^iR$ and
$\widehat{S}:=\lim\limits_{\leftarrow i} S/z^iS$ be the corresponding completions.
The action on $S$ determines an action of $G$ also on $\hS$.
Denote the primitive idempotents of $\hS$ by $e_1,...,e_k$.
For every $1\leq j\leq k$, let $G_j:=\{g\in G| g(e_j)=e_j\}$ be the decomposition group which corresponds to the
primitive idempotent $e_j$.
Let $\widehat{T}:=\lim\limits_{\leftarrow i} T/z^iT$. Then $\hT\simeq\hS*G$,
where the action and 2-cocycle in the crossed product of
$G$ over $\hS$ are induced from those in $T$. The crossed product $\hT$ is an $\hR $-order.
For every $1\leq j\leq k$, let $T_j:=e_j\hT e_j$. Then by the definition of the decomposition groups,
$T_j=\hS e_j*G_j$.
Since $\hS e_j$ is a DVR, we know how to determine if $T_j$ is a maximal $\hR $-order.
The reduction is established by passing from $T$ to $T_1$ by the following
\begin{Theorem}
The following are equivalent
\begin{enumerate}
\item $T$ is a maximal $R$-order.
\item $\hT$ is a maximal $\hR $-order.
\item $T_1$ is a maximal $\hR $-order.
\end{enumerate}
Furthermore, the number of two sided maximal ideals of the above three algebras is equal.
\end{Theorem}
\begin{proof}
(1)$\Leftrightarrow$(2) The number of maximal two sided ideals of $T$ does not change when passing to the completion
$\hT$. Now, by Theorem \ref{HAG}, it remains to show that $T$ is hereditary if and only if so is $\hT$.
Indeed, since $S/\Jac(S)$ and $\hS/\Jac(\hS)$ are both isomorphic to $k$ copies of the residue field of $S$
(where $\Jac$ denotes the Jacobson radical), we obtain that $\bar{T}:=T/\Jac(S)T\simeq\hT/\Jac(\hS)\hT$.
By \cite[Theorem A]{BGL} (a general version of Theorem \ref{heredity}), both $T$ and $\hT$ are hereditary
if and only if $\bar{T}$ is semisimple. \\
(2)$\Leftrightarrow$(3) Since the action of $G$ on the set $\{e_j\}_{j=1}^k$ of primitive idempotents of $\hS$
is transitive, it follows that for every $1\leq j\leq k$, $e_j\in \hT e_1\hT$. Consequently, $1\in \hT e_1\hT$
and hence $\hT =\hT e_1\hT$. By \cite[Proposition 3.5.6]{MR},
we deduce that $\hT$ and $e_1\hT e_1$ are Morita equivalent and we are done.
\end{proof}


\end{document}